\documentclass[a4paper, oneside, 11pt]{article}
\usepackage[english]{babel}
\usepackage{amsfonts}
\usepackage{amssymb}
\usepackage{graphics}
\usepackage{amsrefs}
\usepackage{latexsym}
\begin{document}

\title{{\bf{\Large{Stirling's approximations for exchangeable Gibbs weights}}}\footnote{{\it AMS (2000) subject classification}. Primary: 60G58. Secondary: 60G09.} }
\author{\textsc {Annalisa Cerquetti}\footnote{Corresponding author, SAPIENZA University of Rome, Via del Castro Laurenziano, 9, 00161 Rome, Italy. E-mail: {\tt annalisa.cerquetti@gmail.com}}\\
\it{\small Department of Methods and Models for Economics, Territory and Finance}\\
  \it{\small Sapienza University of Rome, Italy }}
\newtheorem{teo}{Theorem}
\date{\today}
\maketitle{}

\begin{abstract}
We obtain some approximation results for the weights appearing in the exchangeable partition probability function identifying Gibbs partition models of parameter $\alpha \in (0,1)$, as introduced in Gnedin and Pitman (2006). We rely on approximation results for central and non-central generalized Stirling numbers and on known results for conditional and unconditional $\alpha$ diversity. We provide an application to an approximate Bayesian nonparametric estimation of discovery probability in species sampling problems under normalized inverse Gaussian priors. 
\end{abstract}

\section {Introduction}

Exchangeable Gibbs partitions  (Gnedin and Pitman, 2006) are the largest class of infinite exchangeable partitions of the positive integers $\mathbb{N}$ and are characterized by a consistent sequence of exchangeable partition probability functions (EPPF) in the form 
$$
p(n_1, \dots, n_k)= V_{n,k} \prod_{j=1}^{k} (1- \alpha)_{n_j-1},
$$
for each $n \geq 1$, $(n_1, \dots, n_k)$ the sizes of the blocks, $\alpha \in (-\infty,1)$ and $V_{n,k}$ coefficients satisfying the backward recursive relation
$$
V_{n,k}=(n-k\alpha)V_{n+1, k}+ V_{n+1,k+1},
$$ 
for $V_{1,1}=1$. The specific form of the $V_{n,k}$ identifies the specific $\alpha$-Gibbs model, each arising (cfr. Th. 12, Gnedin and Pitman, 2006) as a mixture of extreme partitions, namely: Fisher's (1943) $(\alpha, \xi)$ partitions, $\xi=1, 2, \dots$, for $\alpha < 0$, Ewens' (1972) Poisson-Dirichlet $(\theta)$ partitions for $\alpha=0$, and Pitman's (2003) Poisson-Kingman $PK(\rho_\alpha|t)$ conditional partitions for $\alpha \in (0,1)$, where $\rho_\alpha$ is the L\'evy density of the Stable subordinator. For $\alpha \in (0,1)$, $\theta >-\alpha$ and mixing density 
$$
\gamma_{\theta, \alpha}(t)= \frac{\Gamma(\theta +1)}{\Gamma(\theta/\alpha +1)} t^{-\theta} f_\alpha(t),
$$  
where $f_\alpha(t)$ is the density of the stable distribution, the Poisson-Kingman $(\rho_\alpha, \gamma_{\alpha, \theta})$ model corresponds to the {\it two-parameter Poisson-Dirichlet} partition model (Pitman, 1995, Pitman and Yor, 1997) with EPPF in the form
\begin{equation}
\label{alfateta}
p_{\alpha, \theta}(n_1, \dots, n_k)= \frac{(\theta + \alpha)_{k-1 \uparrow \alpha}}{(\theta +1)_{n-1}} \prod_{j=1}^{k} (1- \alpha)_{n_j-1},
\end{equation}
which reduces to the Dirichlet $(\theta)$ partition model for $\alpha=0$. Recently exchangeable Gibbs partitions and the corresponding discrete random probability measures, have found application  as prior models in the Bayesian nonparametric treatment of species sampling problems, (cfr. Lijoi {\it et al.} 2007b, Lijoi {\it et al.} 2008).  Here interest typically lies in estimating the diversity of a population of species with unknown relative abundances, by estimating both predictive {\it species richness} (cfr. Favaro et al. 2009) and posterior {\it species evenness} (cfr. Cerquetti, 2012). By assuming the sequence of labels of different observed species to be a realization of an exchangeable sequence $(X_i)_{i\geq 1}$ with de Finetti measure belonging to the Gibbs class, posterior predictive results are  obtained with respect to a future $m$-sample conditioning on the multiplicities  $(n_1, \dots, n_k)$  of the first $k$ species observed in a basic $n$-sample. \\\\
With respect to the whole $PK(\alpha, \gamma)$ class,  the $(\alpha, \theta)$ partition model (\ref{alfateta}) stands out for its mathematical tractability, and a huge amount of results are available for this model both in random partitions theory like in Bayesian nonparametric implementations. But in general, without loss of generality, for mixing distributions in the form $\gamma (t)=h(t)f_{\alpha}(t)$, for $f_\alpha(t)$ the Stable density, the Gibbs weights are obtained by mixing the conditional EPPF $p_\alpha(n_1, \dots, n_k|t)$, (cfr Pitman, 2003, eq. (66)) by the specific mixing distribution, which results in the following integral
\begin{equation}
\label{pesigen}
V_{n,k}^{\alpha, h} = \int_0^{\infty}\frac{\alpha^k}{\Gamma(n -k\alpha)} t^{-k \alpha} \left[\int_0^1 p^{n-1-k\alpha} f_\alpha((1-p)t)dp \right]h(t)dt.
\end{equation}
It follows that calculating explicitly the $PK(\rho_\alpha, h \times f_\alpha)$ coefficients outside the $(\alpha, \theta)$ class can be hard, mostly due to the lack of explicit expression for the $\alpha$-stable density. To give an example, in the Bayesian nonparametric implementation of those models, the most studied alternative to the $(\alpha, \theta)$ model is the {\it normalized generalized Gamma} class (Pitman, 2003; see Lijoi {\it et al.} 2005, 2007a, Cerquetti, 2007) whose mixing density  arises by the exponentialy tilting of the stable density,  hence
$$
\gamma(t)= h(t)\times f_\alpha(t)= \exp\{\psi_\alpha(\lambda) - \lambda t\}f_\alpha(t),
$$
for $\psi_\alpha(\lambda)= (2\lambda)^{\alpha}$, $\lambda >0$ the Laplace exponent of $f_\alpha$.  By the reparametrization  $\lambda= \beta^{1/\alpha}/2$, 
\begin{equation}
\label{ngg_1}
\gamma_{\alpha, \beta}(s^{-1/\alpha})= \exp \left\{\beta -\frac 12 \left(\frac {\beta}{s}\right)^{1/\alpha}\right\}f_\alpha(s^{-1/\alpha}) \alpha^{-1} s^{-1/\alpha -1}
\end{equation}
and the corresponding $(V_{n,k})$' (Cerquetti, 2007) are given by
%\begin{equation}
%\label{weingg}
$$V_{n,k}^{\alpha, \beta}= \frac{e^{\beta} 2^n \alpha^{k}}{\Gamma(n)} \int_{0}^{\infty} \lambda^{n-1} \frac{e^{-(\beta^{1/\alpha} +2\lambda)^{\alpha}}}{(\beta^{1/\alpha} +2\lambda)^{n-k\alpha}} d\lambda,
$$
%\end{equation}
and can be rewritten (Lijoi {\it et al.} 2007a) by the change of variable $x=(\beta^{1/\alpha} + 2\lambda)^{\alpha}$, $d\lambda= (2\alpha)^{-1} x^{1/\alpha -1}dx$, as a linear combination of incomplete Gamma functions
\begin{equation}
\label{ngg_2}
V_{n,k}^{\alpha, \beta}=\frac{e^{\beta} \alpha^{k-1}}{\Gamma(n)} \sum_{i=0}^{n-1} {n -1 \choose i} (-1)^{i}(\beta)^{i/\alpha} \Gamma(k - \frac{i}{\alpha}; \beta).
\end{equation}
The computational burden implicit in working with those priors is then clear. 
%Some results for additionally explicit Gibbs partitions generated by the stable subordinator beyond the generalized Gamma class are in Ho {\it et al.} (2007).  
Additionally, when dealing with posterior predictive analysis in this setting, Bayesian nonparametric distributions for quantities of interest in species sampling problems are typically characterized by a ratio, that we term {\it posterior Gibbs weights}, of the kind
\begin{equation}
\label{posteriorwe}
\frac{V_{n+m, k+k^*}}{V_{n,k}},
\end{equation}
for $k^*$ the number of different species observed in the additional $m$-sample, whose complexity increases with the complexity of the {\it prior} weights. It looks therefore interesting to investigate the possibility to obtain some sort of approximation result for both prior and posterior Gibbs coefficients. Here we propose a Stirling's approximation  by relying on approximation results for both generalized {\it central} and {\it non central} Stirling numbers.
The paper is organized as follows. In Section 2. we provide some preliminaries and obtain an approximation result for the {\it prior} Gibbs weights of the $PK(\rho_\alpha, h \times f_\alpha)$ class. In Section 3 we obtain the approximation result for the posterior Gibbs weights (\ref{posteriorwe}) by relying on a recent result for conditional $\alpha$ diversity (Cerquetti, 2001) and in Section  4. we apply our findings to an approximate posterior estimation of discovery probability under normalized inverse gaussian partition model (Pitman, 2003, Lijoi et al. 2005) which corresponds to the generalized Gamma case for $\alpha=1/2$.

\section {Stirling's  approximations for $PK(\rho_\alpha, h\times f_\alpha)$ weights}

From the first order Stirling's approximation for ratio of Gamma functions  $\frac{\Gamma(n +a)}{\Gamma(n+b)}\sim n^{a-b}$ is an easy task to verify that, for $n \rightarrow \infty$ and $k \approx sn^{\alpha}$, the following approximation holds for the EPPF (\ref{alfateta}) of the $PD(\alpha, \theta)$ model
\begin{equation}
\label{stipd}
p_{\alpha, \theta}(n_1, \dots, n_k)\approx \frac{(\alpha)^{k-1}\Gamma(k)}{\Gamma(n)} \frac{\Gamma(\theta +1)}{\Gamma(\theta/\alpha +1)} \left(\frac{k}{n^{\alpha}}\right)^{\theta/\alpha} \prod_{j=1}^{k} (1-\alpha)_{n_j-1}.
\end{equation}
Our first aim is to generalize (\ref{stipd}) to the $PK(\rho_\alpha, \gamma)$ class. First recall that {\it generalized Stirling numbers} are combinatorial coefficients appearing in the expansion of rising factorials $(x)_n=x(x+1)\cdots(x+n-1)$ in terms of generalized rising factorials $(x)_{k \uparrow \alpha}=x(x+\alpha)(x+2\alpha)\cdots(x +(k-1)\alpha)$, hence
$$
(x)_n= \sum_{k=1}^{n} S_{n,k}^{-1, -\alpha} (x)_{k \uparrow \alpha}.
$$
An explicit expression for those numbers is given by Toscano's (1939) formula
\begin{equation}
\label{explstir}
S_{n,k}^{-1, -\alpha}= \frac{1}{\alpha^k k!} \sum_{j=1}^{k} (-1)^j {k \choose j}(-j\alpha)_{n},
\end{equation}
additionally (cfr. e.g. Pitman, 2006) they correspond to {\it partial Bell polynomials} of the kind
$$
B_{n,k}(1-\alpha_\bullet)= \sum_{\{A_1, \dots, A_k\} \in \mathcal{P}_{[n]^k}} \prod_{j=1}^k (1-\alpha)_{|A_j| -1}= \frac{n!}{k!} \sum_{(n_1, \dots, n_k)} \prod_{j=1}^k \frac{(1-\alpha)_{n_j -1}}{n_j!}
$$
for $\alpha \in (0,1)$, where the first sum is over all {\it partitions} $(A_1, \dots, A_k)$ of $[n]$ in $k$ blocks, and the second sum in over all {\it compositions} of $n$ into $k$ parts.
%It can be showns they correspond to the $(n,k)th$ partial Bell polynomials for the combinatorial structure $(1 -\alpha)_{\bullet}$, which is  the number of ways to partition $[n]$ into $k$ blocks an to assing to each block the $(1-\alpha)_\bullet$ combinatorial structure
The following result relies on an approximation for generalized Stirling numbers obtained in Pitman (1999) and on a local approximation result for the distribution of a proper normalization of the numbers of blocks $K_n$ of an $\alpha \in (0,1)$ Gibbs partition. \\\\
\noindent{\bf Proposition 1.} {\it Let $V_{n,k}^{\alpha, h}$ be the coefficients in the EPPF of an exchangeable Gibbs partition arising by a general $PK(\rho_\alpha, h \times f_\alpha)$ Poisson-Kingman model, then the following Stirling's approximation holds for $n \rightarrow \infty$  and $k \approx sn^{\alpha}$
$$
V_{n,k}^{\alpha, h} \approx \frac{\alpha^{k-1} \Gamma(k)}{ \Gamma(n)}\hspace{0.1cm} h\left[\left(\frac{k}{n^{\alpha}}\right)^{-1/\alpha}\right].
$$
}\\
{\it Proof.} In Pitman (1999, cfr. eq. (96)) an asymptotic formula for the generalized Stirling numbers $S_{n,k}^{-1, -\alpha}$ for $n \rightarrow \infty$ and $0 <s < \infty$ , with $k\approx s n^{\alpha}$ is derived by known local limit approximations for the number of blocks in a partition generated by a $
PD(\alpha, \alpha)$ model by a stable density, namely
%$$
%\mathbb{P}(K_n=k)=\frac{\alpha^{k} k!}{(\alpha)_n} S_{n,k}^{-1, -\alpha} 
%$$
\begin{equation}
\label{appsti}
S_{n,k}^{-1, -\alpha} \approx \frac{\alpha^{1-k}\Gamma(n)}{\Gamma(k)}g_{\alpha}(s) n^{-\alpha},
\end{equation}
where, $g_\alpha(s)= \alpha^{-1}f_\alpha(s^{-1/\alpha})s^{-1-1/\alpha}$. Now, for $K_n$ the number of blocks in a $PK(\rho_\alpha, h\times f_\alpha)$ partition model (Pitman, 2003) almost surely, for $n \rightarrow \infty$
$$
K_n/n^{\alpha} \rightarrow S_a
$$
for $S_\alpha \sim  h(s^{-1/\alpha})g_{\alpha}(s)$, which implies the following local limit approximation for $k \approx sn^{\alpha}$ holds for the distribution of $K_n$
$$
\mathbb{P}(K_n=k)\approx V_{n,k}^{\alpha, h} S_{n,k}^{-1, -\alpha}\approx h(s^{-1/\alpha})g_{\alpha}(s)n^{-\alpha}.
$$
Thus the result follows by substitution. \hspace{8.2cm} $\square$ \\\\
{\bf Example 2.} The result in Proposition 1. agrees with the Stirling's approximation for the weights of the $(\alpha, \theta)$ model obtained in (\ref{stipd}). It is enough to notice that $PD(\alpha, \theta)=PK(\rho_\alpha, \gamma_{\alpha, \theta})$ for $\gamma_{\theta, \alpha}= h_{\alpha, \theta} \times f_\alpha$ corresponding to the  $\theta$ polynomial tilting of the stable density for
$$
h(t)= \frac{\Gamma(\theta +1)}{\Gamma(\theta/\alpha +1)} t^{-\theta}.
$$
As for the Generalized Gamma model recalled in (\ref{ngg_1}) and (\ref{ngg_2}), the first order Stirling's approximation for the Gibbs weights  will correspond to
$$
V_{n,k}^{\alpha, \beta} \approx \frac{\alpha^{k-1}\Gamma(k)}{\Gamma(n)} \exp\left\{\beta -\frac n2 (\beta/k)^{1/\alpha}\right\}.
$$
\section {Stirling's approximation for posterior $PK(\rho_\alpha, h\times f_\alpha)$ weights}
As recalled in the Introduction, posterior predictive distributions in Bayesian nonparametric estimation in species sampling problems under Gibbs priors are typically characterized by a ratio of Gibbs weights. As an example consider the posterior joint distribution of the random vector $K_m, L_m, S_1, \dots, S_{K_m}$ for $K_m$ the number of new species generated by the additional $m$ sample, $L_{m}$ the number of new observations in new species, and  $S_1, \dots, S_{K_m}$ the vector of the sizes of the new species in {\it exchangeable random order}, namely
%arises by multiplying (\ref{oldenew}) for the suitable combinatorial coefficient namely
$$
\mathbb{P} (K_m= k, L_m=s, {\bf S}_{K_m}= (s_1, \dots, s_{K_m})| {\bf n})=
$$
$$
= \frac{s!}{s_1! \cdots s_{k^*}! k!} \frac{V_{n+m, k+k^*}}{V_{n,k}} {m \choose s} (n -k\alpha)_{m-s} \prod_{i=1}^{k^*} (1 -\alpha)_{s_i -1}$$
%or
%\begin{equation}
%\label{newblock}
%=\frac{m!}{s_1! \cdots s_k! k! m-s!}  \frac{V_{n+m, j+k}}{V_{n,j}} (n -j\alpha)_{m-s} \prod_{i=1}^{k} (1 -\alpha)_{s_i -1} 
%\end{equation}
or the Bayesian nonparametric estimator for the number $K_m$ of new species in the additional sample (Lijoi et al. 2007b)
$$
\mathbb{E} (K_{m}|K_n=k)= \sum_{k^*=0}^m k^* \frac{V_{n+m, k+k^*}}{V_{n,k}} S_{m, k^*}^{-1, -\alpha, -(n-k\alpha)}.
$$
Here $S_{m,k^*}^{-1, -\alpha, -(n -k\alpha)}$ are non-central generalized Stirling numbers, defined by the convolution relation (see Hsu and Shiue, 1998, Eq. (16)) 
\begin{equation}
\label{conv_sti}
S_{m,k^*}^{-1, -\alpha, -(n-k\alpha)}= \sum_{s=k^*}^{m} {m \choose s} (n-k\alpha)_{m-s}S_{s, k^*}^{-1, .\alpha}, 
\end{equation}
whose corresponding Toscano's formula (\ref{explstir}) may be derived by (\ref{conv_sti})  as follows
$$
S_{m,k^*}^{-1, -\alpha, -(n -k\alpha)}= \sum_{s=k^*}^{m} {m \choose s} (n-k\alpha)_{m-s} \frac{1}{\alpha^*}\frac{1}{k^*!} \sum_{j=1}^{k^*} (-1)^{j}{k^* \choose j} (-j \alpha)_{s}= 
$$
\begin{equation}
\label{toscgen}
= \frac{1}{\alpha^{k^*}} \frac{1}{k^*!} \sum_{j=1}^{k^*} (-1)^{j} {k^* \choose j} \sum_{s=k^*}^{m} {m \choose s} (n -k\alpha)_{m-s} (-j \alpha)_{s}=\frac{1}{\alpha^{k^*}} \frac{1}{k^*!} \sum_{j=1}^{k^*} (-1)^{j} {k^* \choose j} (n -(j+k)\alpha)_{m}.
\end{equation}
In order to derive an approximation result for $$\frac{V_{n+m,k+k^*}}{V_{n,k}}$$ we  replicate the approach adopted in the previous Section by first looking for an approximation for {\it non central} generalized Stirling numbers by exploiting their definition in terms of convolution of central generalized Stirling numbers (\ref{conv_sti}).\\\\
{\bf Proposition 3.} {\it The following approximation holds for non central generalized Stirling numbers $S_{m, k^*}^{-1, -\alpha, -(n-k\alpha)}$ for $\alpha \in (0,1)$ and $m$ large
\begin{equation}
\label{appnonce}
S_{m, k^*}^{-1, -\alpha, -(n -k\alpha)}\approx \frac{\alpha^{1- k^*} \Gamma(m)}{\Gamma(k^*) \Gamma(n -k\alpha)} m^{n-k\alpha -\alpha} \int_0^1 (1-p)^{n-k\alpha-1} p^{-\alpha -1} g_\alpha(zp^{-\alpha})dp
\end{equation}
for
$$
g_\alpha(s)=f_\alpha(s^{-1/\alpha}) \alpha^{-1} s^{-1-1/\alpha}
$$
and $z=k^*/m^{\alpha}$.}\\\\
{\it Proof.} By the definition of non-central Stirling numbers
$$S_{m,k^*}^{-1, -\alpha, -(n-k\alpha)}= \sum_{s=k^*}^{m} {m \choose s}(n-k\alpha)_{m-s} S_{s, k^*}^{-1, -\alpha}
$$
and by (\ref{appsti})
$$
S_{m,k^*}^{-1, -\alpha, -(n-k\alpha)}=\sum_{s=k^*}^m {m \choose s} \frac{\Gamma(s)}{\Gamma(k^*)}s^{-\alpha} \alpha^{1-k^*} g_\alpha\left(\frac{k^*}{s^{\alpha}}\right) (n -k\alpha)_{m-s}.
$$
This may be rewritten as
$$
S_{m,k^*}^{-1, -\alpha, -(n-k\alpha)}=\alpha^{1- k^*} \sum_{s=k^*}^{m} \frac{m\Gamma(m)}{\Gamma(s) \Gamma(m-s+1)}\frac{\Gamma(s)}{\Gamma(k^*)} \frac{\Gamma(n-k\alpha +m -s)}{\Gamma(n-k\alpha)} s^{-\alpha -1} g_{\alpha} \left( \frac{k^*}{m^{\alpha}}\frac{m^{\alpha}}{s^{\alpha}} \right).
$$
For $z=k^*/m^{\alpha}$, by first order Stirling approximations for ratio of Gamma functions
$$
S_{m,k^*}^{-1, -\alpha, -(n-k\alpha)} \approx \frac{\alpha^{1-k^*} \Gamma(m) m^{-\alpha + n -k\alpha -1}}{\Gamma(n-k\alpha) {\Gamma(k^*)}}\sum_{s=k^*}^{m} \left(\frac{m-s}{m}\right)^{n-k\alpha-1} g_\alpha(z p^{-\alpha}) \left( \frac{s}{m}\right)^{-\alpha-1}=
$$
and by the change of variable $p=s/m$
$$
=\frac{\alpha^{1-k^*} m^{-\alpha + n -k\alpha} \Gamma(m)}{\Gamma(n-k\alpha) \Gamma(k^*)} \int_0^1 (1-p)^{n-k\alpha -1} g_{\alpha}(zp^{-\alpha})p^{-\alpha-1} dp.
$$\hspace{15.5cm}$\square$ \\\\
Now we are in a position to obtain the Stirling's approximation for the ratio $V_{n+m, k+k^*}/V_{n,k}$ resorting to a general result for {\it conditional $\alpha$ diversity} for exchangeable Gibbs partitions driven by the Stable subordinator, recently obtained in Cerquetti (2011). \\\\
%\begin{theorem}
%\label{teo1}
{\bf Proposition 4.} {\it The following Stirling's approximation holds for the general posterior Gibbs weights for a $PK(\rho_\alpha, h \times f_\alpha)$ partition model for large $m$ and $k^* \approx s m^{\alpha}$:
\begin{equation}
\label{apppost}
\frac{V_{n+m, k+k^*}}{V_{n,k}}\approx \frac{h(s^{-1/\alpha}) s^k \alpha^{k^*} \Gamma(k^*) \Gamma(n) m^{-(n-k\alpha)}} {\mathbb{E}_{n,k, \alpha} (h(S^{-1/\alpha}))\Gamma(m) \Gamma(k)}= \frac{\alpha^{k+k^*-1} h(s^{-1/\alpha}) s^{k} \Gamma(k^*) m^{-(n-k\alpha)}}{V_{n,k} \Gamma(m)}.
\end{equation}
}\\
{\it Proof.}
Let $\Pi$ be a $PK(\rho_\alpha, \gamma)$ partition of $\mathbb{N}$ driven by the stable subordinator for some  $0 < \alpha <1$ and some mixing probability distribution that without loss of generality we assume in the form $\gamma(t)=h(t)f_\alpha(t)$. Fix $n \geq 1$ and a partition $(n_1, \dots, n_k)$ of $n$ with $k$ positive box-sizes, then by a result in Cerquetti (2011) 
for $K_m$ the number of new blocks induced by an additional $m$-sample
$$
\frac{K_m}{m^{\alpha}}|(K_n=k) \stackrel{d}\longrightarrow S_{\alpha, h}^{n,k} 
$$
for $S_{\alpha, h}^{n,k}$ having density
\begin{equation}
\label{conddiv}
f_{n,k}^{h,\alpha}(s)= \frac{h(s^{-1/\alpha}) \tilde{g}_{n,k}^\alpha(s)}{\mathbb{E}_{n,k}^{\alpha}[h(S^{-1/\alpha})]},
\end{equation}
for  
%\begin{equation}
%\label{gnka}
$$\tilde{g}_{n,k}^{\alpha}(s)= \frac{\Gamma(n)}{\Gamma(n-k\alpha) \Gamma(k)} s^{k-1/\alpha -1} \int_0^1 p^{n-1-k\alpha} f_\alpha((1-p)s^{-1/\alpha})dp 
$$
%\end{equation} 
the density of the product $Y_{\alpha, k} \times [W]^\alpha$ where  $Y_{\alpha, k}$  has density 
%\begin{equation}
%\label{polmit}
$$
g_{\alpha, k\alpha}(y)=\frac{\Gamma(k\alpha +1)}{\Gamma(k+1)}y^kg_{\alpha}(y)
$$
%\end{equation}
for $g_\alpha(y)=\alpha^{-1}y^{-1-1/\alpha}f_\alpha(y^{-1/\alpha})$, independently of $W \sim \beta(k\alpha, n-k\alpha)$. Additionally  
$$
\mathbb{E}_{n,k}^{\alpha} [h(S^{-1/\alpha})]= V_{n,k,h} \frac{\alpha^{1-k} \Gamma(n)}{\Gamma(k)}.
$$
Hence the following local approximation for the posterior distribution of $K_m$ (Lijoi et al. 2007) holds
$$
\mathbb{P}_{\alpha, h}(K_m = k^*|K_n=k)= \frac{V_{n+m, k+k^*}}{V_{n,k}} S_{m, k^*}^{-1, -\alpha, -(n -k\alpha)} \approx  \frac{h(s^{-1/\alpha}) \tilde{g}_{n,k}^\alpha(s)}{\mathbb{E}_{n,k}^{\alpha}[h(S^{-1/\alpha})]} m^{-\alpha}.
$$
Substituting (\ref{appnonce}) the result easily follows. \hspace{8.2cm} $\square$\\\\
{\bf Example 5.} [Poisson-Dirichlet ($\alpha, \theta$) model] Direct first order Stirling's approximation for the posterior $PD(\alpha, \theta)$ coefficients provides
$$
\frac{V_{n+m, k+k^*}}{V_{n,k}}= \frac{(\theta +k\alpha)_{k^* \uparrow \alpha}}{(\theta +n)_m}\approx \frac{\alpha^{k^*}((k^*)^{\theta/\alpha +k}) \Gamma(k^*) \Gamma(\theta +n)}{\Gamma(m) \Gamma(\theta/\alpha +k) m^{\theta +n}},
$$
while an application of (\ref{apppost}) for $h(s^{-1/\alpha})= \frac{\Gamma(\theta +1)}{\Gamma(\theta/\alpha +1)}s^{\theta/\alpha}$ yields 
$$\frac{V_{n+m, k+k^*}}{V_{n,k}} \approx \frac{\alpha^{k^* +k -1} \left[ \left(\frac{k^*}{m^{\alpha}}\right)^{-1/\alpha}\right]^{-\theta} \left(\frac{k^*}{m^\alpha} \right)^{k} \Gamma(k^*) m^{-(n -k\alpha)}\Gamma(\theta+1)}{ \Gamma(m) \Gamma(\theta/\alpha +1)}
\frac{\Gamma(\theta/\alpha +1) \Gamma(\theta +n)}{\alpha^{k-1}\Gamma(\theta/\alpha +k) \Gamma(\theta +1)}
$$
and the two results agree. \\\\
{\bf Example 6.} [Generalized Gamma $(\beta, \alpha)$ model] The exact form of the ratio in this case is given by
$$
\frac{V_{n+m, k+k^+}}{V_{n,k}}= \frac{\alpha^{k^*}\sum_{i=0}^{n+m-1} {n+m-1 \choose i}(-1)^{i} (\beta)^{i/\alpha} \Gamma(k+k^*-i/\alpha; \beta)}{(n)_m\sum_{i=0}^{n-1} {n-1 \choose i}(-1)^{i} (\beta)^{i/\alpha} \Gamma(k-i/\alpha; \beta)}
$$  
while an application of (\ref{apppost}) with $$h[(k^*/m^{\alpha})^{-1/\alpha}]= \exp\left\{\beta - \frac{m}{2}\left( \frac{\beta}{k^*}\right)^{1/\alpha}\right\}$$ provides
$$
\frac{V_{n+m, k+k^*}}{V_{n,k}}\approx \frac{\alpha^{k^*} \exp\left\{- \frac{m}{2}\left(\frac{\beta}{k^*}\right)^{1/\alpha}\right\} 
\left(\frac{k^*}{m^{\alpha}}\right)^k \Gamma(k^*) m^{-(n -k\alpha)}\Gamma(n)}{\Gamma(m) \sum_{i=0}^{n-1} {n -1 \choose i} (-1)^{i}(\beta)^{i/\alpha} \Gamma(k - \frac{i}{\alpha}; \beta)}.
$$
\section{Approximate estimation of discovery probability under inverse Gaussian partition model} In Lijoi {\it et al.} (2007b) a Bayesian nonparametric estimate, under squared loss function of the probability of observing a new species at the $(n+m+1)th$ draw conditionally on a basic sample with observed multiplicities  $(n_1, \dots, n_k)$, without observing the intermediate $m$ observations has been obtained as 
\begin{equation}
\label{disco}
\hat{D}_m^{n,k} = \sum_{k^*=0}^m \frac{V_{n+m+1, k+k^*+1}}{V_{n,k}} S_{m, k^*}^{-1, -\alpha, -(n -k\alpha)}.
\end{equation}
The authors even derive explicit formulas under Dirichlet, two-parameter Poisson-Dirichlet and normalized Inverse Gaussian priors, which correspond to the generalized Gamma prior model for $\alpha=1/2$. Notice that the larger is the size of the additional sample the more computationally hard would be to calculate explicitly (\ref{disco}). 

By an application of (\ref{apppost}),  an approximate evaluation of (\ref{disco}) for large $m$ would corresponds to  
$$
\hat{D}_m^{n,k} \approx \sum_{k^*=0}^{m} \frac{\alpha^{k^*+k} h(s^{-1/\alpha})s^{k} \Gamma(k^*+1) (m+1)^{-(n-k\alpha)}}{V_{n,k} \Gamma(m+1)} S_{m, k^*}^{-1, -\alpha, -(n -k\alpha)} = 
$$
$$=\frac{\alpha^k  (m+1)^{-(n -k\alpha)}}{V_{n,k} \Gamma(m+1)} \sum_{k^*=0}^{m} h((k^*/m^{\alpha})^{-1/\alpha}) (k^*/m^{\alpha})^{k} \alpha^{k^*}\Gamma(k^*+1) S_{m, k^*}^{-1, -\alpha, -(n-k\alpha)},
$$
which specializes under inverse Gaussian model as follows
$$
\hat{D}_m^{n,k}= \frac{(m+1)^{-(n -k/2)}\Gamma(n)}{\Gamma(m+1)\sum_{i=0}^{n-1} {{n-1} \choose i} (-1)^{i} \beta^{2i} \Gamma(k- 2i; \beta)}\times
$$
$$
\times \sum_{k^*=0}^m \exp\left\{-\frac{m+1}{2}\left( \frac{\beta}{k^*}\right)^{2}\right\} \frac{k^*}{\sqrt{m+1}}S_{m,k^*}^{-1, -1/2, -(n-2k)},
$$
and by the Toscano's formula for generalized Stirling numbers recalled in (\ref{toscgen}) simplifies to
$$
\hat{D}_m^{n,k}= \frac{(m+1)^{-(n -k/2)}\Gamma(n)}{\Gamma(m+1)\sum_{i=0}^{n-1} {{n-1} \choose i} (-1)^{i} \beta^{2i} \Gamma(k- 2i; \beta)}\times
$$
$$
\times \sum_{k^*=0}^m \exp\left\{-\frac{m+1}{2}\left( \frac{\beta}{k^*}\right)^{2}\right\} \frac{k^*}{\sqrt{m+1}} \sum_{j=1}^{k^*} (-1)^{j} {k^* \choose j} (n -(j+k)/2)_{m}.
$$
To have an idea of the reduction of the computational burden obtained provided by the approximation here we report the exact  formula obtained in Lijoi et al. (2007b)
$$
\hat{D}_m^{n,k}= \frac{(-\beta^2)^{m+1}}{(n)_{m+1}} \sum_{k^*=0}^m \frac{\sum_{i=0}^{n+m} {{n+m} \choose i} (-\beta^2)^{-i} \Gamma(k+k^* +1 +2i -2(m+n); \beta)}{\sum_{i=0}^{n-1} {{n-1} \choose i} (-\beta^2)^{-i} \Gamma(k+2 +2i -2n; \beta)}\times
$$
$$
\times \sum_{s=k}^m {m \choose s} {{2s -k -1} \choose {s-1}} \frac{2^{k-2s}\Gamma(s)}{\Gamma(k^*)} {n -j/2}_{m-s}.
$$

\section*{References}
\newcommand{\bibu}{\item \hskip-1.0cm}
\begin{list}{\ }{\setlength\leftmargin{1.0cm}}

\bibu \textsc{Cerquetti, A.} (2007) A note on Bayesian nonparametric priors derived from exponentially tilted Poisson-Kingman models. {\it Stat \& Prob Letters}, 77, 18, 1705--1711.

%\bibu \textsc {Cerquetti, A.} (2009) A Generalized sequential construction of exchangeable Gibbs partitions with application. {\it Proceedings of S.Co. 2009, September 14-16, Milano, Italy. }

%\bibu \textsc {Cerquetti, A.} (2011a) A decomposition approach to Bayesian nonparametric estimation under two-parameter Poisson-Dirichlet priors. {\it Proceedings of ASMDA 2011 - Applied Stochastic Models and Data Analysis. Rome, June, 7-10 2011}

\bibu \textsc {Cerquetti, A.} (2011) Conditional $\alpha$-diversity for exchangeable Gibbs partitions driven by the stable subordinator. {\it Proceeding of the 7th Conference on Statistical Computation and Complex Systems, Padova, Italy, 2011}

\bibu \textsc{Cerquetti, A.} (2012) A Bayesian nonparametric estimator of Simpson's index of evenness under Gibbs priors.	{\it arXiv:1203.1666v1 [math.ST]}

%\bibu \textsc {Chao, A., Bunge, J.} (2002) Estimating the number of species in a stochastic abundance model. {\it Biometrics}, 58, 531--539.

%\bibu \textsc {Cerquetti, A.} (2011) Reparametrizing the two-parameter Gnedin-Fisher partition model in a Bayesian perspective. {\it Proceeding of ISI Dublin, 2011}

%\bibu \textsc {Eberhardt, L. L.} (1969) Some aspects of species diversity models. {\it Ecology}, 50, 3.

%\bibu  \textsc {Charalambides, C. A.} (2005) {\it Combinatorial Methods in Discrete Distributions}. Wiley, Hoboken NJ.

\bibu \textsc {Ewens, W.} (1972) The sampling theory of selectively neutral alleles. {\it Theor. Pop. Biol.}, 3, 87--112.

%\bibu \textsc{Favaro, S., Lijoi, A., Mena, R.H. and Pr\"unster, I.} (2009) Bayesian non-parametric inference for species variety with a two-parameter Poisson-Dirichlet process prior. {\it J. Roy. Stat. Soc. B}, 71, 993-1008.

%\bibu \textsc{Favaro, S., Lijoi, A. and Pr\"unster, I.} (2011) Asymptotics for a Bayesian nonparametric estimator of species variety. {\it Bernoulli} (to appear).

%\bibu \textsc {Favaro, S., Lijoi, A. and Pr\"unster, I.} (2012) Conditional formulae for Gibbs-type exchangeable random partitions. {\it Ann. Appl. Probab.} (to appear)

\bibu \textsc{Fisher, R.A., Corbet, A.S. and Williams, C. B.} (1943) The relation between the number of species and the number of individuals in a random sample of an animal population. {\it J. Animal Ecol.} 12, 42--58.

%\bibu \textsc{Gill, C.A. and Joanes, D. N.} (1979) Bayesian estimation of Shannon's index of diversity. Biometrika, 66, 1, 81-85.

%\bibu \textsc{Ginebra, J.and Puig, X.} (2010) On the measure and the estimation of evenness and diversity. {\it Computational Statistics and Data Analysis}, 54, 2187--2201. 

%\bibu \textsc{Gnedin, A., Haulk, S. and Pitman, J.} (2009) {Characterizations of exchangeable partitions and random discrete distributions by deletion properties}. {\sf http://arxiv.org/abs/0909.3642}

\bibu \textsc{Gnedin, A. and Pitman, J. } (2006) {Exchangeable Gibbs partitions  and Stirling triangles.} {\it J. Math.  Sci.}, 138, 3, 5674--5685. 

%\bibu \textsc{Hansen, B. and Pitman, J.} (2000) Prediction rules for exchangeable sequences related to species sampling. {\it Statistics \& Probability Letters}, 46, 251--256.

%\bibu \textsc{Ho, M-W, James, L.F. and Lau, J.W.} (2008) Explicit Gibbs Chinese Restaurant Process priors. {\it Unpublished manuscript}

%\bibu \textsc{Ho, M-W, James, L.F. and Lau, J.W.} (2007) Gibbs partitions (EPPF's) derived from a stable subordinator are Fox H - And Meijer G - Transforms. arXiv:0708.0619v2 [math.PR]

\bibu \textsc{Hsu, L. C, \& Shiue, P. J.} (1998) A unified approach to generalized Stirling numbers. {\it Adv. Appl. Math.}, 20, 366-384.

%\bibu \textsc{Johnson, N.L. \& Kotz, S.} (1977) {\it Urn models and their application}. Wiley \& Sons.

%\bibu \textsc{Johnson, N. L. \& Kotz, S.} (2005) {\it Univariate Discrete Distributions}. 3rd Edition, Wiley \& Sons.

\bibu \textsc{Lijoi, A. Mena, R.H. and Pr\"unster, I.} (2005) Hierarchical mixture modeling with normalized inverse Gaussian priors. {\it J. Am. Stat. Assoc.}, 100, 1278-1291. 

\bibu \textsc{Lijoi, A., Mena, R.H. and Pr\"unster, I.} (2007a) Controlling the reinforcement in Bayesian non-parametric mixture models. {\it J. Roy, Stat. Soc. B}, 69, 4, 715--740. 

\bibu \textsc{Lijoi, A., Mena, R.H. and Pr\"unster, I.} (2007b) Bayesian nonparametric estimation of the probability of discovering new species.  {\it Biometrika}, 94, 769--786.

\bibu \textsc{Lijoi, A., Pr\"unster, I. and Walker, S.G.} (2008) Bayesian nonparametric estimator derived from conditional Gibbs structures. {\it Ann. Appl. Probab.}, 18, 1519--1547.

%\bibu \textsc {Normand, J.M.} (2004) Calculation of some determinants using the $s$-shifted factorial. {\it J. Phys. A: Math. Gen.} 37, 5737-5762.

%\bibu \textsc{Perman, M., Pitman, J, \& Yor, M.} (1992) Size-biased sampling of Poisson point processes and excursions. {\it Probab. Th. Rel. Fields}, 92, 21--39.

\bibu \textsc{Pitman, J.} (1995) Exchangeable and partially exchangeable random partitions. {\it Probab. Th. Rel. Fields}, 102: 145-158.

\bibu \textsc{Pitman, J.} (1996) Some developments of the Blackwell-MacQueen urn scheme. In T.S. Ferguson, Shapley L.S., and MacQueen J.B., editors, {\it Statistics, Probability and Game Theory}, volume 30 of {\it IMS Lecture Notes-Monograph Series}, pages 245--267. Institute of Mathematical Statistics, Hayward, CA.

\bibu \textsc{Pitman, J.} (1999) Brownian motion, bridge, excursion and meander characterized by sampling at independent uniform times. {\it Electron. J. Probab.}, 4, 1-33.

%\bibu \textsc{Pitman, J.} (1996b) Notes on the two parameter generalization of Ewens random partition structure. {\it Manuscript} University of California, Berkeley. Unpublished.

\bibu \textsc{Pitman, J.} (2003) {Poisson-Kingman partitions}. In D.R. Goldstein, editor, {\it Science and Statistics: A Festschrift for Terry Speed}, volume 40 of Lecture Notes-Monograph Series, pages 1--34. Institute of Mathematical Statistics, Hayward, California.

\bibu \textsc{Pitman, J.} (2006) {\it Combinatorial Stochastic Processes}. Ecole d'Et\'e de Probabilit\'e de Saint-Flour XXXII - 2002. Lecture Notes in Mathematics N. 1875, Springer.

\bibu \textsc{Pitman, J. and Yor, M.} (1997) The two-parameter Poisson-Dirichlet distribution derived from a stable subordinator. {\it Ann. Probab.}, 25, 855--900.

%\bibu \textsc {Simpson, E.H.} (1949) Measurement of diversity. {Nature} 163, 688

\bibu \textsc{Toscano, L.} (1939) Numeri di Stirling generalizzati operatori differenziali e polinomi ipergeometrici. {\it Comm. Pontificia Academica Scient. } 3:721-757.

%\bibu \textsc{Yamato, H. and Sibuya, M.} (2000) Moments of some statistics of Pitman sampling formula. {\it Bull. Inform. Cybernet.}, 32 1--10.

\end{list}

\end{document}